\documentclass{amsart}
\usepackage{amsmath, amsthm, amscd, amsfonts, amssymb, graphicx, color}
\usepackage[bookmarksnumbered, colorlinks, plainpages]{hyperref}
\newtheorem{theorem}{Theorem}[section]

\newtheorem{proposition}[theorem]{Proposition}

\newtheorem{corollary}[theorem]{Corollary}
\newtheorem{definition}[theorem]{Definition}

\newtheorem{remark}[theorem]{Remark}

\newcommand{\M}{\mathbb{M}}

\usepackage{comment}
\excludecomment{mysection}

\begin{document}
	
\title{isometries of absolute order unit spaces}
\author{Anil Kumar Karn and Amit kumar}
	
\address{School of Mathematical Sciences, National Institute of Science Education and Research, HBNI, Bhubaneswar, P.O. - Jatni, District - Khurda, Odisha - 752050, India.}

\email{\textcolor[rgb]{0.00,0.00,0.84}{anilkarn@niser.ac.in, amit.kumar@niser.ac.in,}}
	
	
	
	
\begin{abstract}
	We prove that for a bijective, unital, linear map between absolute order unit spaces is an isometry if, and only if, it is absolute value preserving. We deduce that, on (unital) $JB$-algebras, such maps are precisely Jordan isomorphisms. Next, we introduce the notions of absolutely matrix ordered spaces and absolute matrix order unit spaces and prove that for a bijective, unital, linear map between absolute matrix order unit spaces is a complete isometry if, and only if, it is completely absolute value preserving. We obtain that on (unital) C$^*$-algebras such maps are precisely C$^*$-algebra isomorphism. 
\end{abstract}

\thanks{The second author was financially supported by the Senior Research Fellowship of the University Grants' Commission of India.}

\maketitle 

\section{Introduction}
In 1941, Kakutani proved that an abstract $M$-space is precisely a concrete $C(K, \mathbb{R})$ space for a suitable compact and Hausdorff space $K$ \cite{Kak}. In 1943, Gelfand and Naimark proved that an abstract (unital) commutative C$^*$-algebra is precisely a concrete $C(K, \mathbb{C})$ space for a suitable compact and Hausdorff space $K$ \cite{GN}. Thus Gelfand-Naimark theorem for commutative C$^*$-algebras, in the light of Kakutani theorem, yields that the self-adjoint part of a commutative C$^*$-algebra is, in particular, a vector lattice. 
	
On the other hand, Kadison's anti-lattice theorem suggest that the self-adjoint part of a general C$^*$-algebra can not be a vector lattice \cite{Kad51}. Nevertheless, the order structure of a C$^*$-algebra has many other properties which encourages us to expect a `non-commutative vector lattice' or a `near lattice' structure in it. Keeping this point of view, the first author introduced the notion of absolutely ordered spaces and that of an absolute order unit spaces \cite{K18}. The self-adjoint parts of unital C$^*$-algebras and (unital) $M$-spaces are examples of absolute order unit spaces. It was shown that under an additional condition (see \cite[Theorem 4.12]{K16}) an absolutely ordered space turns out to be a vector lattice. One can easily show that under the same condition, an absolute order unit space becomes an $M$-space. Therefore, an absolutely ordered space may be termed as a `non-commutative vector lattice'.
	
For an element $a$ in a C$^*$-algebra $A$, we define the `absolute value' of $a$ as $\vert a \vert := (a^* a)^{\frac{1}{2}}$ and for an element $v$ in a vector lattice $V$, we define the `absolute value' of $v$ as $\vert v \vert := v \vee (-v)$. We recall that for a pair of positive elements $a$ and $b$ in $A$, we have $a b = 0$ if, and only if, $|a - b| = a + b$. Also, for a pair of positive elements $u$ and $v$ in $V$, we have $u \wedge v = 0$ if, and only if, $\vert u - v \vert = u + v$. Thus in both the cases, we can say that $a \perp b$ if, and only if, $\vert a - b \vert = a + b$. In other words, the two kinds of orthogonality relate to the same kind of relation in terms of absolute value. The definition of an absolutely ordered space is influenced by some of the basic properties of the absolute value which hold in the both kind of above-mentioned ordered spaces. 
	
In \cite{K51}, Kadison characterized bijective linear isometries between unital C$^*$-algebras. Since then, many generalizations and extensions of this result has been studied. On the one side, surjective isometries of C$^*$-algebras have been characterized as Jordan triple preserving maps, or these results have been extended up to $JB^*$-triples. In another direction, Jordan isomorphisms have been characterized in terms of absolute value preserving maps together with some or the other conditions. (See, for example, \cite{CW04, G79, K83,MY78, RST01}.) A matricial version of the results of this type were studied by Blecher et al. in \cite{BH02, BL02}.
	
In this paper, we study absolute value preserving maps between absolute order unit spaces. We prove that for a bijective, unital, linear map between absolute order unit spaces is an isometry if, and only if, it is absolute value preserving (Theorem \ref{15}). We deduce that on (unital) $JB$-algebras such maps are precisely Jordan isomorphisms (Corollary \ref{22}). Besides this, we study some elementary properties of absolute value preserving maps. Next, we introduce the notions of absolutely matrix ordered spaces and absolute matrix order unit spaces in the context of matrix ordered spaces and present a matricial version of these results. We prove that for a bijective, unital, $\ast$-linear map between absolute matrix order unit spaces is a complete isometry if, and only if, it is completely absolute value preserving (Theorem \ref{25}). From here, we prove that on (unital) C$^*$-algebras such maps are precisely C$^*$-algebra isomorphism (Corollary \ref{26}). This result was proved in \cite[Corollary 3.2]{BH02}. (Also see, \cite{BL02}.) We give a simple, order-theoretic proof using a trick which is apparently new.

\section{Absolute value preservers on absolutely ordered spaces}

We begin by recalling some basic order theoretic notions. Let $V$ be a real vector space. A non-empty subset $V^+$ of $V$ is called a cone if $V^+$ is closed under vectors' addition as well as scalar multiplication with non-negative real numbers. In this case, $(V, V^+)$ is called a \emph{real ordered vector space}. Also, then $(V, \leq)$ is a partially ordered space with the partial order $u \leq v$ if $v - u \in V^+$ in a unique way, in the sense that $u + w \leq v + w$ and $k u \leq k v$ whenever $u \leq v$, $w \in V$ and $k$ is a positive real number. The cone $V^+$ is said to be \emph{proper}, if $V^+ \cap - V^+ = \{ 0 \}$. It is said to be \emph{generating}, if $V = V^+ - V^+$. A positive element $e \in V^+$ is said to be an order unit for $V$ if for each $v \in V$, there is a positive real number $k$ such that $k e \pm v \in V^+$. The cone $V^+$ is said to be \emph{Archimedean}, if for any $v \in V$ with $k u + v \in V^+$ for a fixed $u \in V^+$ and all positive real numbers $k$, we have $v \in V^+$. 

Let $W$ be a vector subspace of $V.$ Then $W$ is said to be an order ideal of $(V, V^+)$ if, whenever $v\in V^+, w \in W$ with $v \leq w,$ we have $v \in W.$ 

Let $(V, V^+)$ be a real ordered vector space with an order unit $e$ such that $V^+$ is Archimedean. Then $e$ determines a norm on $V$ given by 
$$\| v \| := \inf \{ k > 0: k e \pm v \in V^+ \}$$
in such a way that $V^+$ is norm-closed and for each $v \in V$, we have $\| v \| e \pm v \in V^+$. In this case, we say that $V$ is an \emph{order unit space} and denote it by $(V, e)$.

Now, we recall the notion of absolutely ordered spaces which was introduced by the first author as a possible non-commutative model for vector lattices.

\begin{definition} \cite[Definition 3.4]{K18}
Let $(V, V^+)$ be a real ordered vector space and let $\vert\cdot\vert: V \to V^+$ be a mapping satisfying the following conditions:               
     \begin{enumerate}
          \item[(a)] $\vert v \vert = v$ if $v \in V^+;$
          \item[(b)] $\vert v \vert \pm v \in V^+$ for all $v \in V;$
          \item[(c)] $\vert k \cdot v \vert = \vert k \vert \cdot \vert v \vert$ for all $v \in V$ and $k \in \mathbb{R};$
          \item[(d)] If $u, v$ and $w \in V$ with $\vert u - v \vert = u + v$ and $0 \leq w \leq v,$ then $\vert u - w \vert = u + w;$
          \item[(e)] If $u, v$ and $w \in V$ with $\vert u - v \vert = u + v$ and $\vert u - w \vert = u + w,$ then $\vert u - \vert v \pm w \vert \vert = u + \vert v \pm w \vert.$ 
     \end{enumerate}                 
Then $(V, V^+, \vert\cdot\vert)$ is said to be an \emph{absolutely ordered space}. 
\end{definition} 

\begin{definition}
Let $(V, V^+, \vert\cdot\vert)$ be an absolutely ordered space. Let $W$ be a vector subspace of $V$ and put $W^+ := \lbrace v \in V^+ : v \in W \rbrace$. Then $W$ is said to be an absolutely ordered subspace of $(V, V^+, \vert\cdot\vert)$ if $\vert w \vert \in W^+$ for all $w \in W.$ A vector subspace $W$ of $V$ which is an order ideal of $(V, V^+)$ and an absolutely ordered subspace of $(V, V^+, \vert\cdot\vert)$ is called absolutely order ideal of $(V, V^+, \vert\cdot\vert)$.
\end{definition}

\begin{remark}\label{0}
		Let $(V, V^+, \vert\cdot\vert)$ be an absolutely ordered space.
	\begin{enumerate}
		\item  The cone $V^+$ is proper and generating. In fact, if $\pm v \in V^+$, then by (a) and (c), we get 
		$$v = \vert v \vert = \vert - v \vert = - v$$
		so that $v = 0$. Next, by (b), for any $v \in V$, we have $$v = \frac{1}{2}\left( (\vert v \vert + v) - (\vert v \vert - v) \right) \in V^+ - V^+.$$ 
		
		\item Let $u, v \in V$ be such that $\vert u - v\vert = u + v$. Then $u, v \in V^+$. For such a pair $u, v \in V^+$, we shall say that $u$ is \emph{orthogonal} to $v$ and denote it by $u \perp v$.
		\item We write, $v^+ := \frac{1}{2}(\vert v \vert + v)$ and $v^- := \frac{1}{2}(\vert v \vert - v)$. Then $v^+ \perp v^-, v = v^+ - v^-$ and $\vert v \vert = v^+ + v^-$. This decomposition is unique in the following sense: If $v = v_1 - v_2$ with $v_1 \perp v_2$, then $v_1 = v^+$ and $v_2 = v^-$.  In other words, every element in $V$ has a unique orthogonal decomposition in $V^+$.
		\item Let $W$ be a vector subspace of $V.$ Then $W$ is absolutely order ideal of $(V, V^+, \vert\cdot\vert)$ if and only if $w^+, w^- \in W^+$ for all $w \in W.$
	\end{enumerate}
\end{remark}

\begin{definition}
	Let $V$ and $W$ be absolutely ordered spaces. A linear map $\phi: V \to W$ is said to be an \emph{absolute value preserving map} ($\vert\cdot\vert$-preserving map, in short), if $\vert \phi(v) \vert = \phi(\vert v \vert)$ for all $v \in V.$ 
\end{definition}
The next result is a tool to use $\vert\cdot\vert$-preserving maps. 

\begin{proposition}\label{16}
	Let $V$ and $W$ be absolutely ordered spaces and let $\phi: V \to W$ be a linear map. Then the following statements are equivalent:
	\begin{enumerate}
		\item[(1)]$\phi$ is $\vert\cdot\vert$-preserving;
		\item[(2)]$\phi \geq 0$ and $\phi(v_1) \perp \phi(v_2)$ for all $v_1, v_2 \in V^+$ with $v_1 \perp v_2$;
		\item[(3)]$\phi \geq 0$ and $\phi(v^+) = \phi(v)^+$ for all $v \in V$;
		\item[(4)]$\phi \geq 0$ and $\phi(v^-) = \phi(v)^-$ for all $v \in V$.
	\end{enumerate}
	\begin{proof}
		(1)$\implies$(2): Let $v_1 \perp v_2.$ Put $v = v_1 - v_2.$ Then $\vert v \vert = v_1 + v_2.$ Since $\phi$ is an additive $\vert\cdot\vert$-map (by (1)), we get $\phi(v_1) + \phi(v_2) = \phi(\vert v \vert) = \vert \phi(v) \vert = \vert \phi(v_1) - \phi(v_2) \vert.$ Thus $\phi(v_1), \phi(v_2) \in W^+$ with $\phi(v_1) \perp \phi(v_2).$
		
		(2)$\implies$(3): Let $v \in V.$ Then $v^+ \perp v^-$ so that by (2), $\phi(v^+) \perp \phi(v^-).$ As $\phi(v) = \phi(v^+) - \phi(v^-),$ we get $\phi(v)^+ = \phi(v^+), \phi(v)^- = \phi(v^-).$
		
		(3)$\implies$(4): If we use the fact, $v^- = (- v)^+.$
		
		(4)$\implies$(1): Let $v \in V.$ Then $\vert v \vert = v^+ + v^- = (- v)^- + v^-.$ Thus by (4), we get  
		$$\phi(\vert v \vert) = \phi(- v)^- + \phi(v)^- = \phi(v)^+ + \phi(v)^- = \vert \phi(v) \vert.$$
	\end{proof}
\end{proposition}

\begin{theorem}\label{12}
	Let $V$ and $W$ be absolutely ordered spaces and let $\phi: V \to W$ be a linear $\vert\cdot\vert$-preserving map. Then 
	\begin{enumerate}
		\item[(1)]$\ker(\phi)$ is an absolutely order ideal of $V$.
		\item[(2)] $\phi(V)$ is an absolutely ordered subspace of $W.$ In particular, $\phi(V)^+ = \phi(V^+)$. 
		\item[(3)] For each $v \in V$, we define $\vert v+ \ker(\phi) \vert = \vert v \vert + \ker(\phi).$ Then $$(V\big/ \ker(\phi), \left(V\big/ \ker(\phi) \right)^+, \vert\cdot\vert)$$ 
		is also an absolutely ordered space, where  
		$$\left(V\big/ \ker(\phi) \right)^+ := \lbrace v + \ker(\phi) : v \in V^+ \rbrace.$$
	\end{enumerate} 
	\begin{proof}
		\begin{enumerate}
			\item[(1)]Let $v \in \ker(\phi).$ Then $\phi(v) = 0$ so that $0 = \vert \phi(v) \vert = \phi(\vert v \vert).$ Thus $\vert v \vert \in \ker(\phi)$ and consequently, $\ker(\phi)$ is an absolutely ordered subspace of $V.$ Now, as $\phi \geq 0,$ it follows that $\ker(\phi)$ is an order ideal. 
			\item[(2)]Let $w \in \phi(V),$ say $w = \phi(v)$ for some $v \in V.$ Then $\phi(v^+) = \phi(v)^+ = w^+$ and $\phi(v^-) = \phi(v)^- = w^- \in \phi(V).$ Thus $\phi(V)$ is an absolutely ordered subspace of $W.$ Next, if $w \in \phi(V)^+,$ then $\phi(v)^- = w^- = 0$ so that $w = \phi(v) = \phi(v^+).$ Thus $\phi(V)^+ \subset \phi(V^+).$ Now, being $\vert\cdot\vert$-preserving, $\phi \geq 0$ so that $\phi(V^+) \subset \phi(V)^+$. Hence $\phi(V^+) = \phi(V)^+.$   
			\item[(3)] By \cite[Proposition II.1.1]{MA71}, we know that 
			$\left(V\big/ \ker(\phi) \right)^+$ is a proper cone of $V \big/ \ker(\phi).$
		\end{enumerate} 
		\begin{enumerate}
			\item[(a)]Let $v \in V$ with $v + \ker(\phi) \in \left(V \big/ \ker(\phi) \right)^+.$ There exists $v_0 \in V^+$ such that $v + \ker(\phi) = v_0 + \ker(\phi).$ Thus  $$\vert v + \ker(\phi) \vert = \vert v_0 + \ker(\phi) \vert = \vert v_0 \vert + \ker(\phi) = v_0 + \ker(\phi) = v + \ker(\phi).$$
			\item[(b)]Let $v \in V.$ Then 
			$$\vert v + \ker(\phi) \vert \pm (v + \ker(\phi)) = (\vert v \vert \pm v) + \ker(\phi) \in \left(V \big/ \ker(\phi) \right)^+.$$
			\item[(c)]Let $k \in \mathbb{R}.$ Then 
			\begin{eqnarray*}
			\vert k(v + \ker(\phi)) \vert &=& \vert (kv + \ker(\phi)) \vert = \vert k v \vert + \ker(\phi)  \\ 
			&=& \vert k  \vert \vert v \vert + \ker(\phi) = \vert k \vert (\vert v \vert + \ker(\phi)) \\
			&=& \vert k \vert \vert v + \ker(\phi) \vert.
			\end{eqnarray*}
			\item[(d)]Let $u, v, w \in V$ such that $u + \ker(\phi), v + \ker(\phi), w + \ker(\phi) \in \left( V \big/ \ker(\phi) \right)^+$ with $\vert u - v \vert + \ker(\phi) = u + v + \ker(\phi)$ and $w + \ker(\phi) \leq v + \ker(\phi).$ Then $\vert \phi(u - v) \vert = \phi(\vert u - v \vert) = \phi(u + v)$ and $0 \leq \phi(w) \leq \phi(v)$. Since $\phi(V)$ is an absolutely ordered space, we may conclude that 
			$$\phi(\vert u - w \vert) = \vert \phi(u) - \phi(w) \vert = \phi(u) + \phi(w).$$ 
			Thus $\vert u - w \vert + \ker(\phi) = u + w + \ker(\phi)$. 
	 		\item[(e)] Let $\vert u - v \vert + \ker(\phi) = u + v + \ker(\phi)$ and $\vert u - w \vert + \ker(\phi) = u + w + \ker(\phi)$. Then $\vert \phi(u) - \phi(v) \vert = \phi(u) + \phi(v)$ and $\vert \phi(u) - \phi(w) \vert = \phi(u) + \phi(w)$. Since $\phi(V)$ is an absolutely ordered space, we may conclude that $\vert \phi(u) - \vert \phi(v) \pm \phi(w) \vert \vert = \phi(u) + \vert \phi(v) \pm \phi(w) \vert$. Thus, it follows that $\vert u - \vert v \pm w \vert \vert + \ker(\phi) = u + \vert v \pm w \vert + \ker(\phi)$. 
		
			Hence $(V \big/ \ker(\phi), \left(V \big/ \ker(\phi) \right)^+, \vert\cdot\vert)$ is an absolutely ordered space. 
		\end{enumerate}
	\end{proof} 
\end{theorem}

\begin{corollary}\label{14}
	Let $V$ and $W$ be absolutely ordered spaces and let $\phi: V \to W$ be a linear $\vert\cdot\vert$-preserving map. Put $\ker^+ (\phi) := \lbrace v \in V^+ : \phi(v)= 0 \rbrace,$ then
	\begin{enumerate}
		\item[(1)]$\phi$ is injective if, and only if, $\ker^+ (\phi) = \lbrace 0 \rbrace$
		\item[(2)]$\phi$ is surjective if, and only if, $\phi(V^+) = W^+.$
		\item[(3)]The quotient map $\tilde{\phi}: V \big/ \ker(\phi) \to \phi(V)$ is a bijective $\vert\cdot\vert$-preserving map. 
	\end{enumerate}
	\begin{proof}
		\begin{enumerate}
			\item[(1)]In Theorem \ref{12}(1), we have actually proved that $\ker(\phi) = \ker^+ (\phi) - \ker^+ (\phi).$ Now, the proof is immediate. 
			\item[(2)]If $\phi$ is surjective, it follows, from Theorem \ref{12}(2), that $\phi(V^+) = W^+.$ Conversely, assume that $\phi(V^+) = W^+.$ If $w \in W,$ by assumption there exist $v_1, v_2 \in V^+$ such that $\phi(v_1) = w^+, \phi(v_2) = w^-.$ Put $v = v_1 - v_2$ so that $\phi(v) = w.$ Hence $\phi$ is surjective.
			\item[(3)]This is an immediate consequence of Theorem \ref{12}(3).  
		\end{enumerate} 
	\end{proof}
\end{corollary}
\begin{corollary}
	Let $\phi: V \to W$ be a bijective, linear and $\vert\cdot\vert$-preserving map. Then $\phi^{-1}$ is also a bijective, linear and $\vert\cdot\vert$-preserving map.
	\begin{proof}
		Since $\phi$ is a surjective $\vert\cdot\vert$-preserving map, by Corollary \ref{14}(2), we have $\phi(V^+) = W^+.$ Now, as $\phi(V^+) = W^+$ and $\phi$ is injective, we get that  $\phi^{-1} \geq 0.$  Next, let $w_1, w_2 \in W^+$ with $w_1 \perp w_2.$ Then there exist $v_1, v_2 \in V$ such that $\phi(v_i) = w_i, i = 1, 2.$ Put $v = v_1 - v_2$ and consider $v = v^+ - v^-.$ By Corollary \ref{14}(3) and \ref{14}(4), we get $\phi(v^+) = \phi(v)^+ = (w_1 - w_2)^+ = w_1,$ $\phi(v^-) = \phi(v)^- = w_2.$ By injectivity of $\phi, ~ v_1 = v^+$ and $v_2 = v^-$ so that $\phi^{-1}(w_1) \perp \phi^{-1}(w_2).$ By Proposition \ref{16}(2), we conclude that $\phi^{-1}$ is also an $\vert\cdot\vert$-preserving map.
	\end{proof} 
\end{corollary}

\section{Absolute value preservers on absolute order unit spaces}
We begin this section by recalling the notion of an absolute order unit space. First, we consider three types of orthogonality in an absolutely ordered space. 
\begin{definition}[\cite{K18}, Definition 3.6]
	Let $(V, V^+, \vert\cdot\vert)$ be an absolutely ordered space and let $\| \cdot \|$ be a norm on $V$. 
	\begin{enumerate} 
		\item[(a)] For $u, v \in V^+$, we say that $u$ is \emph{orthogonal} to $v$ ($u \perp v$) if, $\vert u - v \vert = u + v$ (Remark \ref{0}(2));
		\item[(b)] For $u, v \in V^+$, we say that $u$ is \emph{$\infty$-orthogonal} to $v$ ($u \perp_\infty v$) if, $\| \alpha u + \beta v \| = max(\| \alpha u \|, \| \beta v \|)$ for all $\alpha, \beta \in \mathbb{R};$
		\item[(c)] For $u, v \in V^+$, we say that $u$ is \emph{absolutely $\infty$-orthogonal} to $v$ ($u \perp_\infty^a v$) if, $u_1 \perp_\infty v_1$ whenever $0 \leq u_1 \leq u$ and $0 \leq v_1 \leq v.$
	\end{enumerate}
\end{definition}
In an absolutely ordered space these three types of orthogonality are related in the following way. 
\begin{proposition}[\cite{K18}, Proposition 3.7]
	Let $(V, V^+, \vert\cdot\vert)$ be an absolutely ordered space and assume that $\Vert\cdot\Vert$ 
	is a norm on $V$ such that $V^+$ is $\Vert\cdot\Vert$-closed. Then the following conditions are equivalent:
	\begin{enumerate}
		\item[(A)] For each $v \in V$, we have 
		$$ \Vert \vert v \vert \Vert = \Vert v \Vert = \max (\Vert v^+ \Vert, \Vert v^- \Vert);$$ 
		\item[(B)] For $u, v \in V^+$, we have $u \perp_{\infty}^a v$ whenever $u \perp v;$
		\item[(C)] For $u, v \in V^+$, we have $u \perp_{\infty} v$ whenever $u \perp v$.
	\end{enumerate}
	If $\Vert\cdot\Vert$ is an order unit norm determined by the order unit $e$, then the above conditions are also equivalent to:
	\begin{enumerate}
		\item[(D)] For each $v \in V$ with $\pm v \leq e$, we have $\vert v \vert \leq e$.
	\end{enumerate}
\end{proposition}
This result leads to the following 
\begin{definition}[\cite{K18}, Definition 3.8]
Let $(V, V^+, |\cdot|)$ be an absolutely ordered space and let $\| \cdot \|$ be an order unit norm on V determined by the order unit $e$ such that $V^+$ is $\| \cdot \|$-closed. Then $(V, V^+, |\cdot|, e)$ is called an \emph{absolute order unit space}, if the following two conditions are satisfied:
\begin{enumerate}
\item[(a)] $\| v \| \leq \max (\|u\|, \|w\|)$ for all $u,v$ and $w \in V$ with $u \leq v \leq w;$
\item[(b)] $\perp = \perp^a_\infty$ on $V^+.$
\end{enumerate}
\end{definition}
Note that the self-adjoint part of a unital C$^*$-algebra is an absolute order unit space. More generally, every unital $JB$-algebra is also an absolute order unit space.

\begin{definition}[\cite{K10}, Definition 3.2]
Let $V$ and $W$ be absolutely ordered spaces with norms $\| \cdot \|_V$ and $\| \cdot \|_W$ respectively and let $\phi: V \to W$ be a linear map. We say that $\phi$ is an order isometry, if $l_W(\phi(v)) = l_V(v)$ for all $v\in V$, where 
$$l_V(v) = \inf \lbrace \| u \|_V: u \in V^+, u + v \in V^+ \rbrace.$$ 
If, in addition, $\phi$ is surjective, we say that $V$ is order isometric to $W$.  
\end{definition}
\begin{proposition}\label{21}
	Let $(V, e_V)$ and $(W, e_W)$ be absolute order unit spaces and let $\phi: V \to W$ be a surjective linear map. Then $\phi$ is an order isometry if, and only if it is a unital isometry.
	\begin{proof}
		First, let $\phi$ be an order isometry. Then, by \cite[Proposition 3.5]{K10}, it is an isometry. Next, put $\phi(e_V) = w_0 \in W^+.$ Then $\| w_0 \| = \| e_V \| = 1$ so that $0 \leq w_0 \leq e_W.$ Put $\phi^{-1}(e_W) = u_0.$ As $\phi(V^+) = W^+$ and $\phi$ is one-one, $u_0 \in V^+.$ Also $\| u_0 \| = \| e_W \| = 1.$ Thus $0 \leq u_0 \leq e_V.$ But then $e_W = \phi(u_0) \leq \phi(e_V) = w_0 \leq e_W$ so that $e_W = w_0 = \phi(e_V)$. Thus $\phi$ is unital.
		
		Conversely, assume that $\phi$ is a unital isometry. We show that it is positive. Let $v \in V^+$ with $\| v \| \leq 1$. Then $0 \leq v \leq e_V$ so that $0 \leq e_V - v \leq e_V$. Thus $\| e_V - v \| \leq 1$. Since $\phi$ is a unital isometry, we get $\| e_W - \phi(v) \| \leq 1$. Then $- e_W \leq e_W - \phi(v) \leq e_W$ so that $\phi(v) \in W^+$. Thus $\phi$ is positive. Since $\phi$ is a unital surjective isometry, we get that $\phi^{-1}$ is also positive. Now, for any $v \in V$, we have  
		\begin{eqnarray*}
			l_V(v) &=& \inf \lbrace \| u \|: u \in V^+, u + v \in V^+ \rbrace \\ 
			&=& \inf \lbrace \| \phi(u) \|: \phi(u) \in W^+, \phi(u) + \phi(v) \in W^+ \rbrace \\
			&\geq& l_W(\phi(v))
		\end{eqnarray*} 
		and dually, 
		\begin{eqnarray*}
			l_W(\phi(v)) &=& \inf \lbrace \| w \|: w \in W^+, w + \phi(v) \in W^+ \rbrace \\ 
			&=& \inf \lbrace \| \phi^{-1}(w) \|: \phi^{-1}(w) \in V^+, \phi^{-1}(w) + v \in V^+ \rbrace \\
			&\geq& l_V(v).
		\end{eqnarray*} 
	Thus $\phi$ is an order isometry.
	\end{proof}
\end{proposition}

\begin{theorem}\label{15}
Let $(V, e_V)$ and $(W, e_W)$ be absolute order unit spaces and let $\phi: V \to W$ be a linear bijection. Then $\phi$ is a unital $\vert\cdot\vert$-preserving map if and only if it is a (surjective) order isometry.
\begin{proof}
First, assume that $\phi$ be a unital $\vert\cdot\vert$-preserving map. Since $\phi$ is surjective, by Corollary \ref{14}(2), $\phi(V^+) = W^+$.  As it is also an injection, it follows that $r e_V \pm v \in V^+$ if and only if $r e_W \pm \phi(v) \in W^+$. Thus $\| v \| = \| \phi(v) \|$ for all $v \in V$ so that by Proposition \ref{21}, $\phi$ is an order isometry. 

Conversely, let $\phi$ be a (surjective) order isometry. Then, by Proposition \ref{21}, $\phi$ is a unital isometry. Next, we show that $\phi$ preserves $\vert \cdot \vert.$ Let $v_1, v_2 \in V^+$ with $v_1 \perp v_2.$ If $v_1 = 0$ or $v_2 = 0$, then $\phi(v_1) \perp \phi(v_2).$ Now, assume that $v_1 \neq 0, v_2 \neq 0.$ Then $w_i = \phi(v_i) \in W^+ \setminus \lbrace 0 \rbrace, i = 1, 2.$ Let $0 \leq u_i \leq w_i, i = 1, 2.$ Then $0 \leq \phi^{-1}(u_i) \leq v_i, i = 1, 2.$ Since, $v_1 \perp v_2$, we have $v_1 \perp_{\infty}^a v_2$ and consequently, $\phi^{-1}(u_1) \perp_{\infty} \phi^{-1}(u_2)$. Thus, by \cite[Theorem 3.3]{K14}, we have 
\begin{eqnarray*}
1 &=& \left\| \| \phi^{-1}(u_1) \|^{-1} \phi^{-1}(u_1) + \| \phi^{-1}(u_2) \|^{-1} \phi^{-1}(u_2) \right\| \\
&=& \left\| \|u_1 \|^{-1} u_1 + \| u_2 \|^{-1} u_2 \right\|.
\end{eqnarray*}
as $\phi^{-1}$ is an isometry. Again, applying \cite[Theorem 3.3]{K14}, we get that $u_1 \perp_{\infty} u_2$ so that $w_1 \perp_{\infty}^a w_2$. Now, by the definition of an absolute order unit space, we get $w_1 \perp w_2$. Hence, by Proposition \ref{16}, $\phi$ is $|\cdot|$-preserving.   
\end{proof}
\end{theorem}
Maitland Wright and Youngson proved that any surjective linear unital isometry $\phi: A \to B$ between unital $JB$-algebras $A$ and $B$ is a Jordan isomorphism \cite[Theorem 4]{MY78}. If we club this result with Theorem \ref{15}, we may deduce the following: 
\begin{corollary}\label{22}
	Let $A$ and $B$ be unital $JB$-algebras and let $\phi: A \to B$ be a bijective linear map. Then the following statements are equivalent:
	\begin{enumerate} 
		\item $\phi$ is an order isometry;
		\item $\phi$ a unital $\vert\cdot\vert$-preserving map;  
		\item $\phi$ is a Jordan isomorphism.
	\end{enumerate}
	\begin{proof}
		It follows from Proposition \ref{21} that (1) implies (3). Next, let $\phi$ be a Jordan isomorphism. Let $\phi(1_A) = p \in B$ and let $\phi^{-1}(1_B) = q \in A$. Then 
		$$1_B = \phi(q) = \phi(1_A o q) = \phi(1_A) o \phi(q) = p o 1_B = p$$ 
		so that $\phi$ is unital.  Also, $\phi$ is positive. In fact, if $a \in A^+$, then $a = (a^{\frac{1}{2}})^2$ so that $\phi(a) = \phi(a^{\frac{1}{2}})^2 \in B^+$. Now, for any $x \in A$, we have 
		$$\vert \phi(x) \vert^2 = \phi(x)^2 = \phi(x^2) = \phi(\vert x \vert^2) = \phi(\vert x \vert)^2$$ 
		so that $\phi(\vert x \vert) = \vert \phi(x) \vert$ for all $x \in A$. Thus (3) implies (2). Now, by Theorem \ref{15}, the proof is complete.
	\end{proof}
\end{corollary}
  
Let $(V, e)$ be an absolute order unit space and let $u, v \in V^+$. We say that $u$ is \emph{absolutely compatible} with $v$ (we write, $u \bigtriangleup v$) if $\vert u - v \vert + \vert e - u - v \vert = e.$. Let us recall the notion of order projection of V, given in \cite[Definition 5.2]{K18}: Let $0 \leq p \leq e$. We say that $p$ is an \emph{order projection}, if $p \perp e - p$. We write $\mathcal{OP}(V)$ for the set of all order projections in $V$. Recall that in a unital C$^*$-algebra, an order projection is precisely a projection \cite[Theorem 5.3]{K18}.
\begin{proposition}
	Let $(V, e_V)$ and $(W, e_W)$ be absolute order unit spaces. Then a unital $\vert\cdot\vert$-preserving map $\phi: V \to W$ preserves order projections.
	\begin{proof}
		Let $p \in \mathcal{OP}(V)$. Then $p \perp e_V - p$. As $\phi(e_V) = e_W$ and $\phi(\vert v \vert) = \vert \phi(v) \vert$ for all $v \in V$, by Proposition \ref{16}, we get $\phi(p) \perp e_W - \phi(p)$. Thus $\phi(p) \in \mathcal{OP}(W)$.
	\end{proof}
\end{proposition}

\begin{theorem}
Let $V$ and $W$ be absolute order unit spaces and let $\phi: V \to W$ is a $\vert\cdot\vert$-preserving map such that $\phi(e_V) \in OP(W)$. Then $\phi(u)$ is absolutely compatible with $\phi(v)$, whenever $u, v \in V^+$ with $u$ is absolutely compatible with $v.$
\begin{proof}
Let $u, v \in V^+$ such that $u \bigtriangleup v.$ Then $\vert u - v \vert + \vert e_V - u - v \vert = e_V.$ It follows, from \cite[Proposition 4.2]{K18}, that $u, v \in [0, e_V]$. As  $\phi$ preserves $\vert\cdot\vert$, it is a positive map so that $\phi(u), \phi(v) \in [0, \phi(e_V)]$. Also, we get   
\begin{equation}\label{i}
\vert \phi(u) - \phi(v) \vert + \vert \phi(e_V) - \phi(u) - \phi(v) \vert = \phi(e_V). 
\end{equation}
Since $\phi(e_V) \in OP(W)$, we have $\phi(e_V) \perp e_W - \phi(e_V)$. As $\phi(u), \phi(v) \in [0, \phi(e_V)]$, we get that $\phi(u) \perp e_W - \phi(e_V)$ and $\phi(v) \perp e_W - \phi(e_V)$. Now, by the definition of an absolutely ordered space, we get $\vert \phi(e_V) - \phi(u) - \phi(v) \vert \perp e_W - \phi(e_V)$. Thus 

\begin{eqnarray*}
\vert e_W - \phi(u) - \phi(v) \vert &=& \vert (e_W - \phi(e_V)) + (\phi(e_V) - \phi(u) - \phi(v)) \vert \\
&=& (e_W - \phi(e_V)) + \vert \phi(e_V) - \phi(u) - \phi(v) \vert.
\end{eqnarray*}
Therefore, by  (\ref{i}), we get
\begin{eqnarray*}
	\vert \phi(u) &-& \phi(v) \vert + \vert e_W - \phi(u) - \phi(v) \vert \\
	&=& \vert \phi(u) - \phi(v) \vert + (e_W - \phi(e_V)) + \vert \phi(e_V) - \phi(u) - \phi(v) \vert \\
	&=& e_W.
\end{eqnarray*}
Hence $\phi(u) \bigtriangleup \phi(v).$
\end{proof}
\end{theorem}

\begin{remark}
Let $V$ and $W$ be absolute order unit spaces and let $\phi: V \to W$ is a surjective order isometry. For $u, v \in [0, e_V]$, we have $\phi(u) \bigtriangleup \phi(v)$ if, and only if, $u \bigtriangleup v.$ 
\end{remark}

\section{A matricial version of absolute value preserving maps}

Let us recall the following notion introduced by Choi and Effros in \cite{CE77}. A \emph{matrix ordered space} is a $*$-vector space $V$ together with a sequence $\lbrace M_n(V)^+ \rbrace$ with $M_n(V)^+ \subset M_n(V)_{sa} (:= \lbrace v \in M_n(V) : v = v^* \rbrace)$ for each $n \in \mathbb{N}$ satisfying the following conditions: 
\begin{enumerate}
	\item[(a)] $(M_n(V)_{sa}, M_n(V)^+)$ is a real ordered vector space, for each $n \in \mathbb{N}$; and  
	\item[(b)] $\alpha^* v \alpha \in M_m(V)^+$ for all $v \in M_n(V)^+$, $\alpha \in \M_{n,m}$ and $n ,m \in \mathbb{N}$. 
\end{enumerate} 
It is denoted by $(V, \lbrace M_n(V)^+ \rbrace)$. If, in addition, $e \in V^+$ is an order unit in $V_{sa}$ such that $V^+$ is proper and $M_n(V)^+$ is Archimedean for all $n \in \mathbb{N}$, then $V$ is called a \emph{matrix order unit space} and is denoted by $(V, \lbrace M_n(V)^+ \rbrace, e)$.

Now, we introduce the following. 
\begin{definition}\label{17}
Let $(V, \lbrace \ M_n(V)^+ \rbrace)$ be a matrix ordered space and assume that  $\vert\cdot\vert_{m,n}: M_{m,n}(V) \to M_n(V)^+$ for $m, n \in \mathbb{N}$. Let us write $\vert\cdot\vert_{n,n} = \vert\cdot\vert_n$ for every $n \in \mathbb{N}$. Then $\left(V, \lbrace M_n(V)^+ \rbrace, \lbrace \vert\cdot\vert_{m,n} \rbrace \right)$ is called an \emph{absolutely matrix ordered space}, if it satisfies the following conditions: 
\begin{enumerate}
\item[$1.$] For all $n \in \mathbb{N}$, $(M_n(V)_{sa}, M_n(V)^+, \vert\cdot\vert_n)$ is an absolutely ordered space;
\item[$2.$] For $v \in M_{m,n}(V), \alpha \in M_{r,m}$ and $\beta \in M_{n,s},$ we have
$$\vert \alpha v \beta \vert_{r,s} \leq \| \alpha \| \vert \vert v \vert_{m,n} \beta \vert_{n,s};$$
\item[$3.$] For $v \in M_{m,n}(V)$ and $w \in M_{r,s}(V),$ we have
$$\vert v \oplus w\vert_{m+r,n+s} = \vert v \vert_{m,n} \oplus \vert w \vert_{r,s}.$$
Here $v \oplus w := \begin{bmatrix} v & 0 \\ 0 & w \end{bmatrix}$.
\end{enumerate} 
\end{definition}

\begin{proposition}\label{20}
Let $(V, \lbrace M_n(V)^+ \rbrace, \lbrace \vert\cdot\vert_{m,n} \rbrace)$ be an absolutely matrix ordered space.  
\begin{enumerate}
\item[$1.$] If $\alpha \in M_{r,m}$ is an isometry i.e. $\alpha^* \alpha = I_m,$ then $\vert \alpha v \vert_{r,n} = \vert v \vert_{m,n}$ for any $v \in M_{m,n}(V).$
\item[$2.$] If $v \in M_{m,n}(V),$ then $\left\vert \begin{bmatrix} 0_m & v \\ v^* & 0_n \end{bmatrix} \right\vert_{m+n} = \vert v^* \vert_{n,m} \oplus \vert v \vert_{m,n}.$
\item[$3.$] $\begin{bmatrix} \vert v^* \vert_{n,m} & v \\ v^* & \vert v \vert_{m,n} \end{bmatrix} \in M_{m+n}(V)^+$ for any $v \in M_{m,n}.$
\item[$4.$] $\vert v \vert_{m,n} = \left\vert \begin{bmatrix} v \\ 0 \end{bmatrix} \right\vert_{m+r,n}$ for any $v \in M_{m,n}$ and $r \in \mathbb{N}.$
\item[$5.$] $\vert v \vert_{m,n} \oplus 0_s = \left\vert \begin{bmatrix} v & 0 \end{bmatrix} \right\vert_{m,n+s}$ for any $v \in M_{m,n}(V)$ and $s \in \mathbb{N}.$
\end{enumerate} 
\begin{proof}
\begin{enumerate}
\item[(1)] Let $\alpha \in M_{r,m}$ be an isometry. Then, using Definition \ref{17}(2), we have 
$$\vert \alpha v \vert_{r,n} \leq \| \alpha \| \vert v \vert_{m,n} = \vert \alpha^* \alpha v \vert_{m,n} \leq \| \alpha^* \| \vert \alpha v \vert_{r,n} = \vert \alpha v \vert_{r,n}.$$ 
Thus $\vert \alpha v \vert_{r,n} = \vert v \vert_{m,n}.$
\item[(2)] Put $\alpha = \begin{bmatrix} 0_{n,m} & I_n \\ I_m & 0_{m,n} \end{bmatrix} \in M_{n+m}.$ Then $\alpha$ is an isometry with $$\alpha \begin{bmatrix} 0_m & v \\ v^* & 0_n \end{bmatrix} = \begin{bmatrix} v^* & 0_n \\ 0_m & v \end{bmatrix}.$$ Now, by (1) and Definition \ref{17}(3), we get $$\left\vert \begin{bmatrix} 0_m & v \\ v^* & 0_n \end{bmatrix} \right\vert_{m+n} = \left\vert \begin{bmatrix} v^* & 0_n \\ 0_m & v \end{bmatrix} \right\vert_{m+n} = \vert v^* \vert_{n,m} \oplus \vert v \vert_{m,n}.$$ 
\item[(3)] As $\begin{bmatrix} 0_m & v \\ v^* & 0_n \end{bmatrix} \in M_{m+n}(V)_{sa}$, by the definition of an absolutely ordered space, we have 
$$\begin{bmatrix} \vert v^* \vert_{n,m} & v \\ v^* & \vert v \vert_{m,n} \end{bmatrix} = \left\vert \begin{bmatrix} 0_m & v \\ v^* & 0_n \end{bmatrix} \right\vert_{m+n} + \begin{bmatrix} 0_m & v \\ v^* & 0_n \end{bmatrix} \in M_{m+n}(V)^+.$$
\item[(4)] For $\alpha = \begin{bmatrix} I_m \\ 0 \end{bmatrix} \in M_{m+r,m}$, we have $\alpha v = \begin{bmatrix} v \\ 0 \end{bmatrix} \in M_{m+r, n}$. Since $\alpha^* \alpha = I_m$, by (1) we get, $\vert v \vert_{m,n} = \left\vert \begin{bmatrix} v \\ 0 \end{bmatrix} \right\vert_{m+r,n}$ if, $v \in M_{m,n}(V)$ and $r \in \mathbb{N}$.
\item[(5)] For $\alpha = \begin{bmatrix} I_m \\ 0 \end{bmatrix} \in M_{m+r,m}$, we have $\alpha \begin{bmatrix} v & 0 \end{bmatrix} = \begin{bmatrix} v & 0 \\ 0 & 0 \end{bmatrix} \in M_{m+r, n+s}$. Since $\alpha^* \alpha = I_m$, by (1), we get  
$$\vert \begin{bmatrix} v & 0 \end{bmatrix} \vert_{m,n+s} = \left\vert \begin{bmatrix} v & 0 \\ 0 & 0 \end{bmatrix} \right\vert_{m+r,n+s} = \begin{bmatrix} \vert v \vert_{m, n} & 0 \\ 0 & 0 \end{bmatrix},$$ 
if $v \in M_{m,n}(V)$ and $r, s \in \mathbb{N}$.
\end{enumerate}
\end{proof}
\end{proposition}

\begin{definition} 
Let $(V, \lbrace M_n(V)^+ \rbrace, e)$ be a matrix order unit space such that 
\begin{enumerate}
\item[(a)] $\left(V, \lbrace M_n(V)^+ \rbrace, \lbrace \vert \cdot \vert_{m,n} \rbrace \right)$ be an absolutely matrix ordered space; and
\item[(b)]$\perp = \perp_{\infty}^a$ on $M_n(V)^+$ for all $n \in \mathbb{N}.$ 
\end{enumerate}
Then $(V, \lbrace M_n(V)^+ \rbrace, \lbrace \vert\cdot\vert_{m,n} \rbrace, e)$ is called an \emph{absolute matrix order unit space}. 
\end{definition}
It is straight forward to check that any unital C$^*$-algebra $(A, 1_A)$ is an absolute matrix order unit space. 
\begin{definition}
Let $\left(V, e_V \right)$ and $\left(W, e_W \right)$ be absolute matrix order unit spaces and let $\phi: V \to W$ be a $\ast$-linear map so that $\phi_n(M_n(V)_{sa}) \subset M_n(W)_{sa}$ for every $n \in \mathbb{N}$. We say that $\phi$ is a \emph{complete order isometry} if, $\phi_n: M_n(V)_{sa} \to M_n(W)_{sa}$ is an order isometry for each $n \in \mathbb{N}$. Further, we say that $\phi$ is a \emph{complete $\vert\cdot\vert$-preserving} if, $\phi_n: M_n(V) \to M_n(W)$ is an $\vert\cdot\vert$-preserving map for each $n \in \mathbb{N}$.
\end{definition}

\begin{theorem}\label{25}
Let $(V, e_V)$ and $(W, e_W)$ be absolute matrix order unit spaces and let $\phi: V \to W$ be a $\ast$-linear surjective isomorphism. Then $\phi$ is a complete order isometry if, and only if, $\phi$ is a unital, complete $\vert\cdot\vert$-preserving map.
\begin{proof}
First, let $\phi$ be a complete order isometry. Fix $n \in \mathbb{N}.$ Then $\phi_n: M_n(V)_{sa} \to M_n(W)_{sa}$ is an surjective order isometry. Thus, by Theorem \ref{15}, we get that $\phi_n(e_V^n) = e_W^n$ and that $\phi_n(\vert v \vert_n) = \vert \phi_n(v) \vert_n$ for all $v \in M_n(V)_{sa}.$ Let $v \in M_n(V).$ Then $\begin{bmatrix} o & v \\ v^* & 0 \end{bmatrix} \in M_{2n}(V)_{sa}$ so that 
$$\phi_{2n} \left(\left\vert \begin{bmatrix} o & v \\ v^* & 0 \end{bmatrix} \right\vert_{2n} \right)=\left\vert \phi_{2n}\left(\begin{bmatrix} o & v \\ v^* & 0 \end{bmatrix} \right) \right\vert_{2n}.$$
Thus, by Proposition \ref{20}(2), we get 
$$\begin{bmatrix} \phi_n(\vert v^*\vert_n) & 0 \\ 0 & \phi_n(\vert v \vert_n) \end{bmatrix} = \begin{bmatrix} \vert \phi(v)^* \vert_n & 0 \\ 0 & \vert \phi(v) \vert_n \end{bmatrix}.$$ 
Therefore, $\phi_n(\vert v\vert_n) = \vert \phi_n(v) \vert_n$ for each $n \in \mathbb{N}$. 

Conversely, assume that $\phi(e_V) = e_W$ and that $\phi$ is a complete $\vert\cdot\vert$-preserving map. Fix $n \in \mathbb{N}$. Then $\phi_n(e_V^n) = e_W^n$ and $\phi_n(M_n(V)_{sa}) \to M_n(W)_{sa}$ is an $\vert\cdot\vert$-preserving map. Now, by Theorem \ref{15}, $\phi_n(M_n(V)_{sa}) \to M_n(W)_{sa}$ is an order isometry for all $n \in \mathbb{N}$. Hence $\phi: V \to W$ is a complete order isometry.
\end{proof}
\end{theorem}
\begin{corollary}\label{26}
	Let $A$ and $B$ be any two unital C$^*$-algebras and let $\phi: A \to B$ be a $\ast$-linear bijective map. Then the following facts are equivalent: 
	\begin{enumerate}
		\item $\phi$ is a complete order isometry;
		\item $\phi$ a unital complete $\vert\cdot\vert$-preserving map;  
		\item $\phi$ is a C$^*$-algebra isomorphism.
	\end{enumerate} 
\begin{proof}
	Following Theorem \ref{25}, it suffices to show that (1) (or, equivalently (2)) implies (3). Let $\phi$ be a complete order isometry. Then $\phi_n: M_n(A)_{sa} \to M_n(B)_{sa}$ is an order isometry for each $n \in \mathbb{N}$. Thus, by Corollary \ref{22}, $\phi_n: M_n(A)_{sa} \to M_n(B)_{sa}$ is a Jordan isomorphism for each $n \in \mathbb{N}$. In particular, $\phi_3(a^2) = \phi_3(a)^2$ for any $a \in M_3(A)_{sa}$. Let $x, y \in A$ and consider $a = \begin{bmatrix} 0 & x & 0 \\ x^* & 0 & y \\ 0 & y^* & 0 \end{bmatrix} \in M_3(A)_{sa}$. Then $\phi_3(a^2) = \phi_3(a)^2$ yields that $\phi(x y) = \phi(x) \phi(y)$. Thus $\phi$ is a C$^*$-algebra isomorphism.
\end{proof}
\end{corollary}
\begin{remark}
	It follows, from Corollary \ref{26}, that a unital surjective $\ast$-linear map between unital C$^*$-algebras is complete isometry if, it is a $3$-isometry.
\end{remark}

\end{document}